\documentclass[letterpaper, 10 pt, conference]{ieeeconf}  

\IEEEoverridecommandlockouts                              
\usepackage{graphicx,%
algorithm,%
algorithmicx,%
algpseudocodex,%
amsmath,%
amsfonts,%
verbatim,%
mathtools,%
enumerate,%
bm,%
hyperref,%
mathabx,%
hhline,%
multirow,%
tikz,%
pgf,%
subcaption,%
braket,%
multirow,%
adjustbox,%
breakcites,%
placeins}

\usepackage{amssymb}
\newtheorem{theorem}{Theorem}
\newtheorem{lemma}[theorem]{Lemma}
\newtheorem{assumption}{Assumption}
\newcommand{\R}{\mathbb{R}}
\newcommand{\pipg}{\textsc{pipg}}
\newcommand{\ecos}{\textsc{ecos}}
\newcommand{\osqp}{\textsc{osqp}}

\newcommand{\mb}[1]{\mathbb{#1}}
\newcommand{\mc}[1]{\mathcal{#1}}

\title{\LARGE \bf
Constraint Preconditioning and Parameter Selection for a First-Order Primal-Dual Method applied to Model Predictive Control
}

\author{Govind M. Chari, Yue Yu, and {Beh\c{c}et}~{A\c{c}\i{}kme\c{s}e}
\thanks{This research was supported by ONR grant N00014-20-1-2288 and Blue Origin, LLC; Government sponsorship is acknowledged.}
\thanks{G. Chari is a PhD student in the Department of Aeronautics \& Astronautics at the University of Washington, Seattle, WA 98195
        {\tt\small gchari@uw.edu}}%
\thanks{Y. Yu is an assistant professor in the Department of Aerospace Engineering and Mechanics at the University of Minnesota Twin Cities, Minneapolis, MN 55455.
        {\tt\small yuey@umn.edu}}%
\thanks{B. {A\c{c}\i{}kme\c{s}e} is a professor in the Department of Aeronautics \& Astronautics at the University of Washington, Seattle, WA 98195
        {\tt\small behcet@uw.edu}}%
}

\begin{document}

\maketitle
\thispagestyle{empty}
\pagestyle{empty}

\begin{abstract}
Many techniques for real-time trajectory optimization and control require the solution of optimization problems at high frequencies. However, ill-conditioning in the optimization problem can significantly reduce the speed of first-order primal-dual optimization algorithms. We introduce a preconditioning technique and step-size heuristic for Proportional-Integral Projected Gradient (PIPG), a first-order primal-dual algorithm. The preconditioning technique, based on the QR factorization, aims to reduce the condition number of the KKT matrix associated with the optimization problem. Our step-size selection heuristic chooses step-sizes to minimize the upper bound on the convergence of the primal-dual gap for the optimization problem. These algorithms are tested on two model predictive control problem examples and show a solve-time reduction of at least 3.6x.
\end{abstract}
\section{INTRODUCTION}

Solving optimal control problems in real-time is at the heart of model predictive control (MPC). This control scheme is optimization-based, where we solve a fixed time horizon optimal control problem recursively to generate state and control trajectories which satisfy relevant system constraints while minimizing a predefined cost function.

In this paper, we consider optimal control problems in the form

\begin{subequations}\label{eq:cannonical-form}
    \begin{align}
        \underset{z}{\mathrm{minimize}} 
        \quad & \frac{1}{2}z^\top Pz + q^\top z \\
        \mathrm{subject \; to} 
        \quad & Hz - g = 0 \label{eq:equality-constraint}\\
        & z \in \mb{D} \label{eq:set-D}
    \end{align}
\end{subequations}

\noindent where $H \in \mathbb{R}^{m \times n}$, $m \leq n$, and $P \succ 0$. Constraint \ref{eq:equality-constraint} contains the linear dynamics and set $\mb{D}$ contains state and control constraints. We assume that set $\mb{D}$ is the Cartesian product of sets that have closed-form projections, such as balls, boxes, second-order cones, half-spaces, and the intersection of a ball and second-order cone \cite{Bauschke2011ConvexAA, bauschke2018projecting}. This assumption applies to many typical state and control constraints in MPC \cite{jerez2014embedded}.

First-order optimization algorithms are popular for solving MPC problems due to their ability to warm-start from any prescribed solution (e.g., one that can be obtained from the previous solution in a recursive MPC scheme), which greatly reduces solve-times  and per-iteration computational cost. Unlike second-order methods, first-order methods do not require expensive matrix factorizations, making them attractive for resource constrained embedded systems \cite{banjac2017embedded}.

However, ill-conditioned optimization problems degrade the performance of methods in this class. Consequently, first-order optimization algorithms use preconditioning and parameter selection to improve their convergence behavior. Preconditioning transforms the given optimization problem into an equivalent problem which is able to be solved in fewer iterations. This is achieved by performing a change of variables or transforming the constraints of the problem to reduce the condition number of a matrix related to the optimization problem, such as the KKT matrix. Whereas preconditioning modifies the given optimization problem, parameter selection modifies the optimization algorithm itself. Parameter selection involves choosing parameters of an optimization algorithm such as step-size or relaxation parameters for fastest convergence \cite{ghadimi2015optimal}.

Modified Ruiz equilibriation \cite{osqp} and the hypersphere preconditioner \cite{abhi2023customized} are two existing preconditioning techniques for first-order primal-dual algorithms. These preconditioners aim to reduce the condition number of the KKT matrix. Both perform a change of variables and additionally, modified Ruiz equilibriation scales the constraints with a diagonal matrix. These two preconditioner transform set $\mb{D}$ in Problem \ref{eq:cannonical-form} since they both employ a change of variables. This requires special structure in the matrices defining  the change of variables in order to preserve closed-form projections onto  $\mb{D}$,  
which 
is needed to solve Problem \ref{eq:cannonical-form} by using Proportional-Integral Projected Gradient ({\pipg}), a first-order primal-dual optimization algorithm \cite{yu2022conic, yu2021mpc, yu2022extrapolated}.

Our first contribution is the development of the \textit{QR preconditioner}. This proposed preconditioner does not perform a change of variables, meaning set $\mb{D}$ is not transformed and consequently closed-form projections are preserved. Instead, the QR preconditioner transforms the equality constraints by scaling them with a triangular matrix to minimize the condition number of the KKT matrix. Since a triangular matrix is used for transforming the constraints as opposed to a diagonal matrix, the QR preconditioner has more degrees of freedom to transform the constraints than modified Ruiz equilibriation. The QR preconditioner requires a matrix factorization to perform this constraint transformation, unlike the aforementioned preconditioners, making it more computationally expensive. However, in our numerical examples the factorization time was a small percentage of the overall solve-time. Additionally, in settings requiring the solution of multiple optimization problems with identical equality constraints, such as quasi-convex optimization and MPC, performing the preconditioning offline and solving the preconditioned problem online will result in improved real-time performance. 

Our second contribution is a step-size selection heuristic for {\pipg}. For fastest convergence, {\pipg} requires carefully selecting step sizes, a process that is currently manually tuned. We eliminate manual tuning of step-sizes by introducing a step-size selection heuristic based on choosing step-sizes to minimize the upper bound on the primal-dual gap.

\subsection{Related work}
In this section we will discuss two existing preconditioning methods.

The hypersphere preconditioner \cite{abhi2023customized} performs a change of variables to make the Hessian of the objective function some positive multiple of the identity matrix, making the objective function perfectly conditioned. This preconditioner avoids transforming the constraints. It is applicable to {\pipg} when the Hessian of the objective function is diagonal or has an easy to compute Cholesky factorization, and the change of variables induced by the preconditioner preserves the ability to project onto set $\mb{D}$ in closed-form.

Modified Ruiz equilibration \cite{osqp} preconditions the objective function and transforms the constraints by performing a linear change of variables with a diagonal matrix and scaling the constraints with a diagonal matrix to make the preconditioned KKT matrix have equal row and column norms, a heuristic for good conditioning. Similar to the hypersphere preconditioner, to be applicable to {\pipg} the change of variables induced must preserve the ability to project onto set $\mb{D}$ in closed-form. 

\subsection{Notation}
We denote the $n \times 1$ vector of zeros as $0_n$, the $m \times n$ matrix of zeros as $0_{m \times n}$, the $n \times n$ identity matrix as $I_n$, the Euclidean projection of point $z$ onto the convex set $\mb{D}$ as $\Pi_{\mb{D}}[z]$, the concatenation of vectors $u$ and $v$ as $(u,v)$, the vector formed from the $i^{th}$ through the $j^{th}$ component of $z$ as $z^{i:j}$, and the Euclidean norm of a vector $z$ as $\|z\| := \sqrt{z^\top z}$.

\section{Proportional-Integral Projected Gradient}\label{sec:pipg}

Proportional-Integral Projected Gradient ({\pipg}), shown in Algorithm \ref{alg:pipg}, is a first-order primal-dual conic optimization algorithm capable of infeasibility detection. This algorithm is specialized to handle common constraints sets which arise in optimal control problems, such as boxes, balls, and second-order cones, via closed-form projections \cite{yu2022conic, yu2021mpc, yu2022extrapolated}. {\pipg}  has been demonstrated to be faster than many state-of-the art optimization solvers such as {\ecos}, \textsc{mosek}, \textsc{gurobi}, and {\osqp} \cite{yu2022extrapolated}. It has also been extensively applied to trajectory optimization problems including three-degree-of-freedom (DoF) powered-descent guidance and has been used as a subproblem solver within sequential convex programming for multi-phase powered-descent guidance, 6-DoF powered-descent-guidance, and spacecraft rendezvous problems \cite{abhi2023customized, elango2022customized, kamath2023seco, chari2024fast, chari2024spacecraft}.

\begin{algorithm}[!ht]
\caption{{\pipg}}
\begin{algorithmic}[1]
\Require \(k_{\max}, \alpha, \beta, z^1\in\mathbb{D}, v^1\).
\Ensure \(z^k\).
\For{\(k=1, 2, \ldots, k_{\max}-1\)}
\State{\(w^{k+1}= v^k+\beta(Hz^k-g)\)} \label{eq:integrator}
\State{\(z^{k+1}=\Pi_{\mathbb{D}}[z^k-\alpha(Pz^k + q +H^\top w^{k+1})]\)} \label{eq:primal-update}
\State{\(v^{k+1}=w^{k+1}+\beta H(z^{k+1}-z^k)\)}
\EndFor
\end{algorithmic}
\label{alg:pipg}
\end{algorithm}

In Algorithm \ref{alg:pipg}, $z^k$ and $w^k$ are the iterates of the primal and dual variables respectively. Thus, we refer to $\alpha$ and $\beta$ as the primal and dual step-sizes. For the convergence results to hold for {\pipg}, the primal and dual step-sizes must satisfy the algebraic relationship

\begin{equation*}
    \alpha(\lambda_{\mathrm{max}} + \sigma \beta) < 1.
\end{equation*}

We parameterize this family of step-sizes with $\gamma > 0$ as

\begin{equation}\label{eq:const-step-size}
    \alpha < \frac{1}{\lambda_{\mathrm{max}} + \gamma}, \quad \beta = \frac{\gamma}{\sigma}
\end{equation}

\noindent where $\lambda_{\mathrm{max}}$ is the maximum eigenvalue of $P$ and $\sigma$ is the maximum eigenvalue of $H^\top H$. In practice we choose $\alpha$ to be equal to the upper bound and do not observe convergence issues. We will discuss how to choose $\gamma$ in Section \ref{sec:parameter-selection}.
\section{QR PRECONDITIONER}\label{sec:qr}
In this section, we derive and provide a geometric interpretation for the QR preconditioner.

\subsection{Derivation}\label{sec:justification}

First, we derive a preconditioner that only modifies Constraint \ref{eq:equality-constraint} of Problem \ref{eq:cannonical-form} and minimizes the condition number of the KKT matrix.
First, we need the following  assumptions on the constraint matrix $H$ and the Hessian of the objective function $P$.

\begin{assumption}\label{assum:rank}
      $H \in \mb{R}^{m \times n}$ has full row rank.  
\end{assumption}

\begin{assumption}\label{assum:P-pd}
    $P$ is positive definite and has maximum and minimum eigenvalues $\lambda_{\mathrm{max}}$ and $\lambda_{\mathrm{min}}$ respectively.
\end{assumption}

The KKT conditions of Problem \ref{eq:cannonical-form} are as follows

\begin{equation}\label{eq:kkt}
    \underbrace{
    \begin{bmatrix}
        P & H^\top \\
        H & 0_{m \times m}
    \end{bmatrix}}_{\mc{K}}
    \begin{bmatrix}
        z^\star \\
        w^\star
    \end{bmatrix} +\begin{bmatrix}
    q\\
    -g
    \end{bmatrix}\in 
    \begin{bmatrix}
         - \mb{N}_{\mb{D}}(z^\star) \\
        \{0_{m}\}
    \end{bmatrix}
\end{equation}

\noindent where $z^\star$ and $w^\star$ are the optimal primal and dual solutions respectively and $\mb{N}_{\mb{D}}(z^\star)$ is the normal cone of set $\mb{D}$ at $z^\star$. We denote the KKT matrix as $\mc{K}$.

If Assumptions \ref{assum:rank} and \ref{assum:P-pd} hold and $\sigma_{\mathrm{max}}$ and $\sigma_{\mathrm{min}}$ are the maximum and minimum singular values of $H$, the eigenvalues of $\mc{K}$ lie in the following set \cite{benzi2005numerical}

\begin{equation*}
    \lambda(\mc{K}) \subset \mc{I}^- \cup \mc{I}^+ \\
\end{equation*}

\noindent where

\begin{equation*}
        \begin{split}
        \mc{I}^- = \left[\frac{1}{2}\left(\lambda_{\mathrm{min}} - \sqrt{\lambda_{\mathrm{min}}^2+4\sigma^2_{\mathrm{max}}}\right), \right. \\
        \left .\frac{1}{2}\left(\lambda_{\mathrm{max}} - \sqrt{\lambda_{\mathrm{max}}^2+4\sigma^2_{\mathrm{min}}}\right)\right] \\
    \end{split}
\end{equation*}

\begin{equation*}
                \mc{I}^+ = \left[\lambda_{\mathrm{min}}, \; \frac{1}{2}\left(\lambda_{\mathrm{max}} + \sqrt{\lambda_{\mathrm{max}}^2+4\sigma^2_{\mathrm{max}}}\right)\right].
\end{equation*}
\noindent The spectral condition number of $\mc{K}$ is defined as \cite{benzi2005numerical}

\begin{equation*}
    \kappa(\mc{K}) = \frac{\mathrm{max} \; |\lambda(\mc{K})|}{\mathrm{min} \; |\lambda(\mc{K})|}.
\end{equation*}

The condition number quantifies the difficulty in terms of the number of iterations to solve Equation \ref{eq:kkt} with an iterative method. Thus, a smaller condition number for $\mc{K}$ results in faster convergence for {\pipg} which solves the KKT conditions in Equation \ref{eq:kkt}.

The upper bound for the $\mc{I}^+$ interval is larger in magnitude than the lower bound for the $\mc{I}^-$ interval since $\lambda_{\mathrm{max}}~\geq~\lambda_{\mathrm{min}} > 0$ and $\sigma_{\mathrm{max}} > 0$ by Assumptions \ref{assum:rank} and \ref{assum:P-pd}. We can construct an upper bound on the condition number of $\mc{K}$ as follows

\begin{equation}\label{eq:cn-bound}
\begin{split}
        \kappa(\mc{K}) \leq \mathrm{max}\left(\frac{\lambda_{\mathrm{max}} + \sqrt{\lambda_{\mathrm{max}}^2+4\sigma^2_{\mathrm{max}}}}{\sqrt{\lambda_{\mathrm{max}}^2+4\sigma^2_{\mathrm{min}}} - \lambda_{\mathrm{max}}}, \right. \\
        \left. \frac{\lambda_{\mathrm{max}} + \sqrt{\lambda_{\mathrm{max}}^2+4\sigma^2_{\mathrm{max}}}}{2 \lambda_{\mathrm{min}}}\right).
\end{split}
\end{equation}

\begin{lemma} \label{lemma:optimal-sig-cn}
    The singular values for $H$ which minimize the condition number of the KKT matrix are $\sigma_{\mathrm{min}} = \sigma_{\mathrm{max}} = \sqrt{\lambda_{\mathrm{max}}\lambda_{\mathrm{min}}+\lambda_{\mathrm{min}}^2}$
\end{lemma}

\begin{proof}
    We will first prove by contradiction that $\sigma_{\mathrm{max}}$ and $\sigma_{\mathrm{min}}$ which minimize the bound on the condition number given by Equation \ref{eq:cn-bound} must be equal.

    Suppose that we have $\sigma_{\mathrm{min}} < \sigma_{\mathrm{max}}$ which minimize the bound on the condition number. Both arguments of the $\mathrm{max}$ operator are monotonically increasing in $\sigma_{\mathrm{max}}$, so we can reduce $\sigma_{\mathrm{max}}$ and further reduce the bound on the condition number which produced a contradiction.

    Since $\sigma_{\mathrm{min}} = \sigma_{\mathrm{max}}$ at optimality, the first argument of the $\mathrm{max}$ operator and its derivative are

    \begin{equation*}
        h(\sigma_{\mathrm{max}}) = \frac{\lambda_{\mathrm{max}} + \sqrt{\lambda_{\mathrm{max}}^2+4\sigma^2_{\mathrm{max}}}}{\sqrt{\lambda_{\mathrm{max}}^2+4\sigma^2_{\mathrm{max}}} - \lambda_{\mathrm{max}}}        
    \end{equation*}

    \begin{equation*}
        h'(\sigma_{\mathrm{max}}) = \frac{-8\lambda_{\mathrm{max}}\sigma_{\mathrm{max}}\sqrt{\lambda_{\mathrm{max}}^2+4\sigma_{\mathrm{max}}^2}}{\left(\sqrt{\lambda_{\mathrm{max}}^2+4\sigma_{\mathrm{max}}^2} - \lambda_{\mathrm{max}}\right)^2}.
    \end{equation*}

    Since $\lambda_{\mathrm{max}} > 0$ and $\sigma_{\mathrm{max}} > 0$, $h'(\sigma_{\mathrm{max}})$ is negative and thus $h(\sigma_{\mathrm{max}})$ is monotonically  decreasing in $\sigma_{\mathrm{max}}$.

    It is also clear that the second argument of the $\mathrm{max}$ operator  in Equation \ref{eq:cn-bound} is monotonically increasing in $\sigma_{\mathrm{max}}$. To minimize the maximum of a monotonically decreasing and monotonically increasing function, we must set the two arguments equal. By setting the arguments equal we obtain the minimum and maximum singular values of $H$ which minimize the KKT condition number

    \begin{equation*}
        \sigma_{\mathrm{max}} = \sigma_{\mathrm{min}} = \sqrt{\lambda_{\mathrm{max}}\lambda_{\mathrm{min}} + \lambda_{\mathrm{min}}^2}.
    \end{equation*}
\end{proof}








If $HH^\top=(\lambda_{\mathrm{min}}\lambda_{\mathrm{max}}+\lambda_{\mathrm{min}}^2)I_m$, then the singular values of $H$ are $\sqrt{\lambda_{\mathrm{max}}\lambda_{\mathrm{min}} + \lambda_{\mathrm{min}}^2}$ which minimize the condition number of the KKT matrix by Lemma \ref{lemma:optimal-sig-cn}. For $HH^\top$ to be some scaling of the identity matrix, it is sufficient for the rows of $H$ to be orthogonal. To do this, we will apply a QR factorization to $H^\top$, and transform Constraint \ref{eq:equality-constraint} into Constraint \ref{eq:preconditioned-constraint} as follows

\begin{subequations}
\begin{align}
   &  Hz-g = 0 \\
   \Rightarrow & \  R^\top Q^\top z - g = 0 \\
   \Rightarrow & \ Q^\top z - R^{-\top}g = 0 \\
   \Rightarrow & \  \eta (Q^\top z - R^{-\top}g) = 0 \label{eq:scale-eta} \\
  \Rightarrow & \  \hat{H}z - \hat{g} = 0 \label{eq:preconditioned-constraint}
\end{align}
\end{subequations}

\noindent where $H^\top = QR$ is the thin or economy QR factorization of $H^\top$, $\hat{H} = \eta Q^\top$, $\hat{g} = \eta R^{-\top}g$, and Assumption \ref{assum:rank} ensures invertibility of $R^\top$ \cite{golub1996matrix}. Note that Assumption \ref{assum:rank} is equivalent to requiring that linear independence constraint qualification (LICQ) of Constraint \ref{eq:equality-constraint} holds for Problem \ref{eq:cannonical-form} \cite{nocedal2006numerical}. We can arbitrarily scale the constraint in Equation \ref{eq:scale-eta} by some $\eta > 0$ and the rows of $\hat{H}$ will still be orthogonal. If we choose $\eta = \sqrt{\lambda_{\mathrm{max}}\lambda_{\mathrm{min}}+\lambda_{\mathrm{min}}^2}$, then $\hat{H}\hat{H}^\top=(\lambda_{\mathrm{min}}\lambda_{\mathrm{max}}+\lambda_{\mathrm{min}}^2)I_m$, thus minimizing the condition number of the KKT matrix. The QR preconditioner is given by Algorithm \ref{alg:qr-preconditioner}.

\begin{algorithm}[!ht]
\caption{QR Preconditioner}
    \begin{algorithmic}[1]
    \Function{qr\_preconditioner}{$H$, $g$, $\lambda_{\mathrm{max}}$, $\lambda_{\mathrm{min}}$}
    \State{\(\eta = \sqrt{\lambda_{\mathrm{max}}\lambda_{\mathrm{min}}+\lambda_{\mathrm{min}}^2}\)}
    \State{\(Q, R = \mathrm{qr}(H^\top)\)} \Comment{economy QR factorization}
    \State{\(\hat{H} = \eta Q^\top\)}
    \State{\(\hat{g} = \eta R^{-\top}g\)}
    \State \Return $\hat{H}$, $\hat{g}$
    \EndFunction
\end{algorithmic}
\label{alg:qr-preconditioner}
\end{algorithm}

\subsection{Geometric Interpretation}

\begin{figure}[!htb]
    \centering
    \includegraphics[width=8cm]{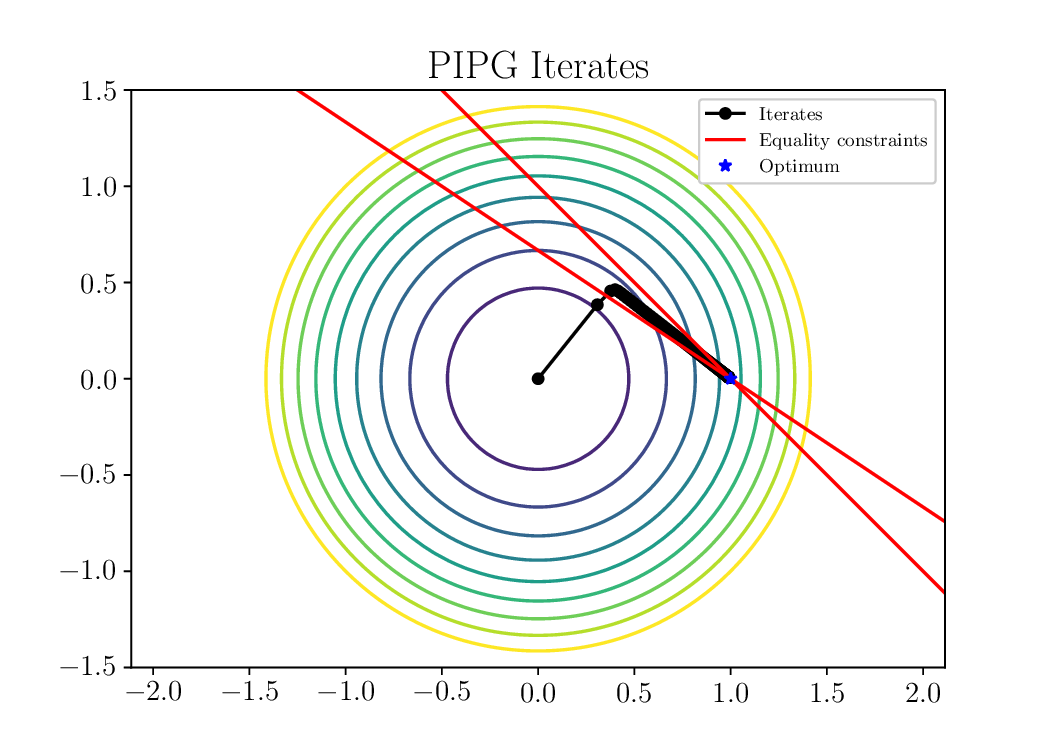}
  
    \caption{{\pipg} iterates with nearly parallel equality constraints}
    \label{fig:qr-motivation}
\end{figure}

\begin{figure}[!htb]
    \centering
    \includegraphics[width=8cm]{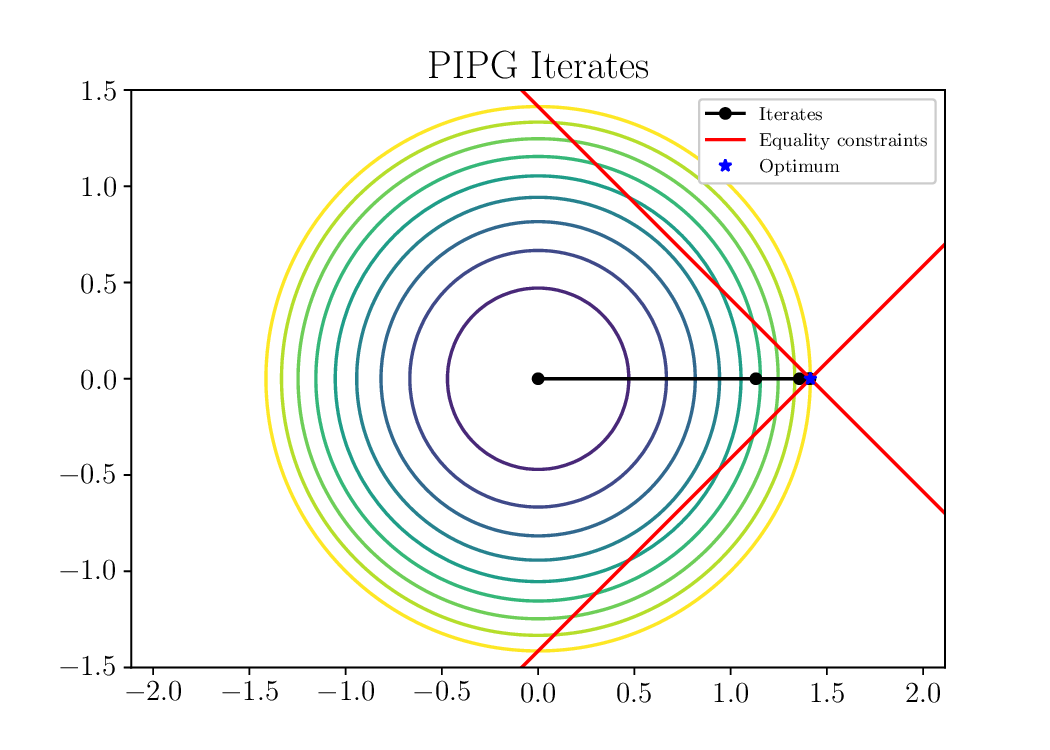}
  
    \caption{{\pipg} iterates after applying QR preconditioner}
    \label{fig:qr-motivation-1}
\end{figure}

To interpret this preconditioner geometrically, we will consider the simplified problem of using {\pipg} to minimize a bivariate quadratic function subject to two equality constraints with nearly parallel normal vectors. Figures \ref{fig:qr-motivation} and \ref{fig:qr-motivation-1} depict the level curves of the objective function and the {\pipg} iterates before and after the QR preconditioner is applied.

We can see that in Figure \ref{fig:qr-motivation}, {\pipg} takes several iterations to converge. Notice that Line \ref{eq:primal-update} of {\pipg} takes a step in the negative gradient direction of the objective function in order to minimize the objective, then takes a step in a direction given by a linear combination of the rows of $H$ (the $H^\top w^{k+1}$ term) to satisfy Constraint \ref{eq:equality-constraint}. The rows of $H$ are the normal vectors of the hyperplanes which define the equality constraints in Problem \ref{eq:cannonical-form}. If we have two equality constraints with nearly parallel hyperplane normals (as in Figure \ref{fig:qr-motivation}), we require smaller coefficients in $w^{k+1}$ for the linear combination of the hyperplane normals, given by $H^\top w^{k+1}$, to move the iterates in a direction nearly parallel to hyperplane normals (up and to the right in Figure \ref{fig:qr-motivation}), but require large coefficients for the iterates to move in a direction nearly orthogonal to the hyperplane normals (down and to the right in Figure \ref{fig:qr-motivation}). The coefficients of this linear combination are given by our current iterate of the dual variables $w^{k+1}$. To get large enough $w^{k+1}$ for the iterates to move nearly orthogonal to the hyperplane normals, we must wait for the integrator given by Line \ref{eq:integrator} to accumulate enough error, resulting in slow convergence.

The QR preconditioner transforms the equality constraints given by Equation \ref{eq:equality-constraint} to a equivalent set of equality constraints with the same feasible set as the original constraints, but one whose hyperplane normal vectors are all orthogonal and of the same length, making it easy for iterates to move in all directions. After applying the QR preconditioner, we can see the resulting iterates in Figure \ref{fig:qr-motivation-1} where the optimum is achieved in far fewer iterations.
\section{PARAMETER SELECTION}\label{sec:parameter-selection}

In this section we will discuss how to select the step-sizes $\alpha$ and $\beta$ to minimize the upper bound on the primal-dual gap for the iterates of {\pipg}.

For convergence of {\pipg}, the primal and dual step-sizes must satisfy the relationship

\begin{equation}\label{eq:alpha-beta-relationship}
    \alpha(\lambda_{\mathrm{max}} + \sigma \beta) < 1 
\end{equation}

\noindent where $\alpha$ and $\beta$ are the primal and dual step-sizes respectively, $\lambda_{\mathrm{max}}$ is the maximum eigenvalue of $P$, and $\sigma$ is the maximum eigenvalue of $H^\top H$.

There are infinite choices for $\alpha$ and $\beta$ which satisfy the above relationship, and it has been observed that for fastest convergence a ``proper" choice for $\alpha$ and $\beta$ must be made. From Proposition 1 of \cite{yu2022conic}, we see that zero primal-dual gap is a sufficient condition for optimality of Problem \ref{eq:cannonical-form}. Thus, we propose to choose step-sizes which minimize the upper bound on the primal-dual gap.

In \cite{yu2022extrapolated}, the authors introduce a positive, scalar parameter, referred to as $\omega$, to select an $\alpha$ and $\beta$ which satisfy the relationship given by Equation \ref{eq:alpha-beta-relationship}. We introduce a different scalar parameter $\gamma$, to choose an $\alpha$ and $\beta$ using Equation \ref{eq:const-step-size}. Using $\gamma$ instead of the previously introduced $\omega$ makes the analysis in this section tractable.

Theorem 1 from \cite{yu2022conic} provides an upper bound on the primal-dual gap for Algorithm \ref{alg:pipg}. We will rewrite this theorem to explicitly show dependence on $\gamma$

\begin{equation}\label{eq:primal-dual-gap}
        \mc{L}(\bar{z}^k, w^\star) - \mc{L}(z^\star, \bar{w}^k) \leq  \frac{f(\gamma)}{k} 
\end{equation}

\begin{equation}\label{eq:primal-dual-gap-constant}
            f(\gamma) = \frac{\lambda_{\mathrm{max}} + \gamma}{2}\|z^1-z^\star\|^2 + \frac{\sigma}{2\gamma}\|v^1-w^\star\|^2 \\
\end{equation}

\begin{equation*}
    \bar{z}^k = \frac{1}{k}\sum_{j=1}^k z^{j+1}, \;            \bar{w}^k = \frac{1}{k}\sum_{j=1}^k w^{j+1}
\end{equation*}

\noindent where $\{z^j, w^j\}_{j=1}^k$ are generated by Algorithm \ref{alg:pipg},  $(z^\star, w^\star)$ is the optimal primal-dual solution of Problem \ref{eq:cannonical-form}, and $\mc{L}(z, w) = \frac{1}{2}z^\top Pz + q^\top z + w^\top (Hz - g)$ is the Lagrangian of Problem \ref{eq:cannonical-form}. Both $\gamma$ and $\sigma$ are defined to be positive real numbers. Under this assumption, we see that $f(\gamma)$ given by Equation \ref{eq:primal-dual-gap-constant} is convex in $\gamma$ and has a minimizer given by 

\begin{equation*}
    \gamma^\star = \sqrt{\sigma}\frac{\|v^1-w^\star\|}{\|z^1-z^\star\|}.
\end{equation*}

Since $\gamma^\star$ is a function of the optimal primal-dual solution $(z^\star, w^\star)$, we cannot directly choose the optimal step-sizes. Instead, we suggest a heuristic where for every fixed number of iterations of {\pipg}, we pick the value for $\gamma$ which minimizes an approximation of Equation \ref{eq:primal-dual-gap-constant}, where the optimal primal-dual solution, $(z^\star, w^\star)$, is replaced with our current iterates of the primal and dual variables, $(z^k, w^k)$. The value of $\gamma$ which minimizes this approximation of Equation \ref{eq:primal-dual-gap-constant} after the $k^{th}$ iteration of {\pipg} is 

\begin{equation*}
    \hat{\gamma}^\star = \sqrt{\sigma}\frac{\|v^1-w^k\|}{\|z^1-z^k\|}.
\end{equation*}

The proposed step-size selection algorithm is given by Algorithm \ref{alg:adaptive-step}. We can combine Algorithms \ref{alg:pipg}, \ref{alg:qr-preconditioner}, and \ref{alg:adaptive-step} into Algorithm \ref{alg:combined} to show how our proposed algorithms should be used within {\pipg}.

Although choosing $\gamma$ to minimize the upper bound on the primal-dual gap given by Equation \ref{eq:primal-dual-gap} only guarantees improving the worst-case performance of {\pipg}, we observe improved convergence in our numerical results.

\begin{algorithm}[!ht]
\caption{Parameter selection for {\pipg}}
    \begin{algorithmic}[1]
    \Function{step\_selection}{$z^1$, $v^1$, $z^k$, $w^k$, $\lambda_{\mathrm{max}}$, $\sigma$}
        \State{\(\hat{\gamma}^\star = \sqrt{\sigma}(\|z^1-z^k\| / \|v^1 - w^k\|)\)}
        \State{\(\alpha = 1 / (\lambda_{\mathrm{max}} + \hat{\gamma}^\star)\)}
        \State{\(\beta = \hat{\gamma}^\star / \sigma\)}
        \State \Return $\alpha$, $\beta$
    \EndFunction
\end{algorithmic}
\label{alg:adaptive-step}
\end{algorithm}

\begin{algorithm}[!ht]
\caption{{\pipg} w/ QR Preconditioner and Parameter Selection}
\begin{algorithmic}[1]
\Require \(k_{\max}, k_{\mathrm{update}}, \alpha_{\mathrm{init}}, \beta_{\mathrm{init}}, z^1\in\mathbb{D}, v^1, H, g\).
\Ensure \(z^k\).
\State{\(\hat{H}, \hat{g} = \textsc{qr\_preconditioner}(H, g, \lambda_{\mathrm{max}}, \lambda_{\mathrm{min}})\)}
\State{\(\alpha, \beta = \alpha_{\mathrm{init}}, \beta_{\mathrm{init}}\)}
\For{\(k=1, 2, \ldots, k_{\max}-1\)}
\If{\(\mathrm{mod}(k, k_{\mathrm{update}}) == 0\)}
    \State{\(\alpha, \beta = \textsc{step\_selection}(z^1, v^1, z^k, w^k, \lambda_{\mathrm{max}}, \sigma)\)}
\EndIf
\State{\(w^{k+1}= v^k+\beta(\hat{H}z^k-\hat{g})\)}
\State{\(z^{k+1}=\Pi_{\mathbb{D}}[z^k-\alpha(Pz^k + q +\hat{H}^\top w^{k+1})]\)}
\State{\(v^{k+1}=w^{k+1}+\beta \hat{H}(z^{k+1}-z^k)\)}
\EndFor
\end{algorithmic}
\label{alg:combined}
\end{algorithm}

\section{NUMERICAL RESULTS}\label{sec:numerical}

In this section, we assess the performance of Algorithms \ref{alg:qr-preconditioner} and \ref{alg:adaptive-step} on two MPC problems: regulation of an oscillating mass-spring system and quadrotor obstacle avoidance, a quadratic program (QP) and a second-order cone program (SOCP) respectively. We first define the two problems, discuss implementation details of our algorithms, and finally present the numerical results.

\subsection{Oscillating Mass Control}
The first MPC problem considers regulation of a one-dimensional mass-spring system to equilibrium by applying forces to each of the masses \cite{wang2010fast}. In this problem there are $N$ masses connected to their neighbors by springs with the masses at the end being connected to a wall with springs.

The \textit{state} at time-step $t$ is given by $x_t = (r_t, v_t)$, where $r_t \in \R^N$ is the displacement of the masses from their respective equilibrium positions, and $v_t \in \R^N$ is the velocity of the masses. The \textit{control} at time-step $t$, $u_t \in \R^N$, is the input forces acting on each of the $N$ masses. We constrain the displacement of the masses (\ref{eq:mass-pos-constraint}), the speed of the masses (\ref{eq:mass-vel-constraint}), and the force applied to each mass (\ref{eq:mass-control-constraint}). We discretize the dynamics with a zero-order-hold and sample-time of $\Delta t$ and can write the problem as follows

\begin{subequations}\label{eq:mass-spring-mpc}
    \begin{align}
        \underset{x_t, u_t}{\mathrm{minimize}} 
        \quad & \frac{1}{2}\left(\sum_{t = 1}^{T} x_t^\top Q_t x_t + \sum_{t = 1}^{T-1}u_t^\top R_t u_t\right) \\
        \mathrm{subject\;to} 
        \quad & x_{t+1} = Ax_t + Bu_t \;  \quad t \in [1, T-1] \label{eq:mass-dynamics} \\
        & x_1 = x_{\mathrm{init}} \\
        & \|r_t\|_\infty \leq r_{\mathrm{max}} \quad\quad\quad\quad\!\!\!   t \in [1, T] \label{eq:mass-pos-constraint} \\
        & \|v_t\|_\infty \leq v_{\mathrm{max}} \quad\quad\quad\quad\!\!\!   t \in [1, T] \label{eq:mass-vel-constraint} \\
        & \|u_t\|_\infty \leq u_{\mathrm{max}} \quad\quad\quad\quad\!\!\!\!   t \in [1, T-1]. \label{eq:mass-control-constraint}
    \end{align}
\end{subequations}

All the masses have unit mass, all the springs have unit spring constant, the following non-dimensional parameters are used,

\begin{equation*}
    \begin{split}
        T &= 30, \; \Delta t = 0.1, \; N = 8, \;   \\
        r_{\mathrm{max}} &= 0.75, \; v_{\mathrm{max}} = 0.75, \; u_{\mathrm{max}} = 0.5 \\
        Q_t &= \mathrm{blkdiag}(I_N, 5I_N), \; R_t = I_N.
    \end{split}
\end{equation*}

Each component of the initial state $x_{\mathrm{init}}$ is drawn from the uniform distribution on the interval [-0.5, 0.5]. 

The dynamics given by Equation \ref{eq:mass-dynamics} can be written as Constraint \ref{eq:equality-constraint} with appropriate choice of $H$ and $g$. Constraints \ref{eq:mass-pos-constraint}, \ref{eq:mass-vel-constraint}, and \ref{eq:mass-control-constraint} can be written as Constraint \ref{eq:set-D}, where $\mb{D}$ is the Cartesian product of boxes.

\subsection{Quadrotor Path Planning}

The second MPC problem considered is the three-degree-of-freedom quadrotor path planning problem with obstacle avoidance. The \textit{state} at time-step $t$ is defined as $x_t = (r_t, v_t)$, where $r_t \in \R^3$ is the position of the quadrotor, and $v_t \in \R^3$ is the velocity of the quadrotor. The \textit{control} at time $t$, $u_t \in \R^3$, is the thrust produced by the quadrotor's propellers.

Unlike the oscillating masses regulation problem, where the goal is to regulate all the masses to their equilibrium positions, the objective in this problem is to track a reference trajectory given by $\hat{x}_t$ which linearly interpolates the initial state $x_{\mathrm{init}}$ and the desired target state $x_{\mathrm{target}}$. However, this reference trajectory intersects the obstacle, forcing the optimizer to compute a trajectory that avoids that obstacle. The dynamics of this system are given by a double integrator in three dimensions subject to gravity. These dynamics are discretized with a zero-order-hold to obtain the discrete-time linear system given by \ref{eq:quad-dynamics}, where 

\vspace{0.25cm}

\begin{equation*}
    A = \begin{bmatrix}
        I_3 & (\Delta t) I_3 \\
        0_3 & I_3
    \end{bmatrix} \hspace{0.25cm}
\end{equation*}

\begin{equation*}
        B = \frac{1}{m}\begin{bmatrix}
        (1/2)(\Delta t)^2 I_3 \\
        (\Delta t) I_3
    \end{bmatrix}
\end{equation*}

\begin{equation*}
        c = \begin{bmatrix}
        0_{2 \times 1} \\
        -(1/2)g(\Delta t)^2 \\
        0_{2 \times 1} \\
        -g(\Delta t)
    \end{bmatrix}
\end{equation*}

\noindent $g = 9.8$ is the non-dimensional gravitational acceleration and $m$ is the mass of the quadrotor.

We impose two control constraints: an upper bound on thrust (\ref{eq:quad-thrust-constraint}) and a thrust pointing constraint (\ref{eq:quad-tilt-constraint}). The thrust pointing constraint is a proxy to constrain the tilt angle of the quadrotor, since we are using a three-degree-of-freedom model and cannot directly constrain the quadrotor's tilt.

A cylindrical keep-out zone can be written as
\begin{equation*}
    \|r_t^{1:2}-r_c\|_2 \geq \rho
\end{equation*}
\noindent where $r_c$ is the position of the center of the obstacle and $\rho$ is the radius of the obstacle. This constraint, however, is nonconvex.

This nonconvex keep-out-zone constraint can be replaced with a convex rotating half-space constraint, as in \cite{yu2021mpc}. This rotating half-space constraint is written as \ref{eq:quad-convex-keepout} with

\begin{equation*}
    \begin{split}
        a_t &= [\cos(\psi t + \phi) \; -\sin(\psi t + \phi)]^\top \\
        b_t &= -a_t^\top r_c - \rho
    \end{split}
\end{equation*}

\noindent where $\psi$ is the rotation rate of the half-space and $\phi$ is the initial phase offset of the half-space. This problem can be written as 

\begin{subequations}
    \begin{align}
        \underset{x_t, u_t}{\mathrm{minimize}}~ 
         & \frac{1}{2}\left(\sum_{t = 1}^{T} (x_t-\hat{x}_t)^\top Q_t (x_t-\hat{x}_t) \right. & & \\
          & \left. \quad + \sum_{t=1}^{T-1} u_t^\top R_t u_t \right) & & \nonumber \\
        \mathrm{subject~to}~& x_{t+1} = Ax_t + Bu_t + c & & \!\!\!\!\!\!\!\!\!\!\!\!\!\!\!t \in [1, T-1] \label{eq:quad-dynamics}\\
        & x_1 = x_{\mathrm{init}} & & \\
        & a_t^\top r_t \leq b_t & & \!\!\!\!\!\!\!\!\!\!\!\!\!\!\!t \in [1, T] \label{eq:quad-convex-keepout}\\
        & \|v_t\|_2 \leq v_{\mathrm{max}} & & \!\!\!\!\!\!\!\!\!\!\!\!\!\!\!t \in [1, T] \label{eq:quad-speed-constraint} \\
        & \|u_t\|_2 \leq u_{\mathrm{max}} & & \!\!\!\!\!\!\!\!\!\!\!\!\!\!\!t \in [1, T-1] \label{eq:quad-thrust-constraint} \\
        & \cos(\theta_{\mathrm{max}})\|u_t\|_2 \leq u_t^\top e & &  \!\!\!\!\!\!\!\!\!\!\!\!\!\!\!t \in [1, T-1] \label{eq:quad-tilt-constraint}
    \end{align}
\end{subequations}

\noindent where $e = [0 \; 0 \; 1]^\top$. Additionally, the following non-dimensional parameters are used

\begin{equation*}
    \begin{split}
        T &= 30, \; \Delta t = 0.2, \; m = 3, \; \psi = -0.5  \\
        \phi &= -\pi/4, \; r_c = [2.5 \; 2.5]^\top, \; \rho = 0.25 \\
        v_{\mathrm{max}} &= 1.5, \; u_{\mathrm{max}} = 35, \; \theta_{\mathrm{max}} = 0.1745 \\
        Q_t &= \mathrm{blkdiag}(2I_3, I_3), \; R_t = 0.5 I_3 \\
        x_{\mathrm{init}} &= [0 \; 0\; 5\; 0\; 0\; 0], \; x_{\mathrm{target}} = [5 \; 5\; 5\; 0\; 0\; 0].
    \end{split}
\end{equation*}

With proper choice for $H$ and $g$, we can rewrite the dynamics given by Constraint \ref{eq:quad-dynamics} in the form of Constraint \ref{eq:equality-constraint}. Constraint \ref{eq:quad-convex-keepout} represents a half-space, Constraint \ref{eq:quad-speed-constraint} represents a ball, and Constraints \ref{eq:quad-thrust-constraint} and \ref{eq:quad-tilt-constraint} represent the intersection of a ball and a second-order cone. We can write all of these constraints as Constraint \ref{eq:set-D}, where $\mb{D}$ is the Cartesian product of the aforementioned sets.

\subsection{Implementation Details}

We generate our numerical results using MATLAB. From Algorithm \ref{alg:qr-preconditioner}, we determine that the maximum eigenvalue of $H^\top H$ is $\eta^2$. If Algorithm \ref{alg:qr-preconditioner} is not used, we use power iteration to compute the maximum eigenvalue of $H^\top H$ \cite{golub1996matrix}. To compute maximum and minimum eigenvalues of $P$, $\lambda_{\mathrm{max}}$ and $\lambda_{\mathrm{min}}$ respectively, we simply take the maximum and minimum diagonal elements of $P$ since $P$ is diagonal in our examples. For general $P \succ 0$, we can compute the maximum eigenvalue with power iteration and minimum eigenvalue with either shifted power iteration or inverse power iteration. In all examples, we use $k_{\mathrm{update}} = 25$ when testing Algorithm \ref{alg:adaptive-step}. When running examples without using Algorithm \ref{alg:adaptive-step}, we set $\gamma = \sigma$ which makes the dual step-size $\beta = 1$ consistent with the implementation in \cite{yu2022conic}. We run the oscillating-mass problem $50$ times with $50$ different initial conditions drawn from the distribution described in the previous section.

\begin{figure}[!htb]
    \centering
    \includegraphics[width=\linewidth]{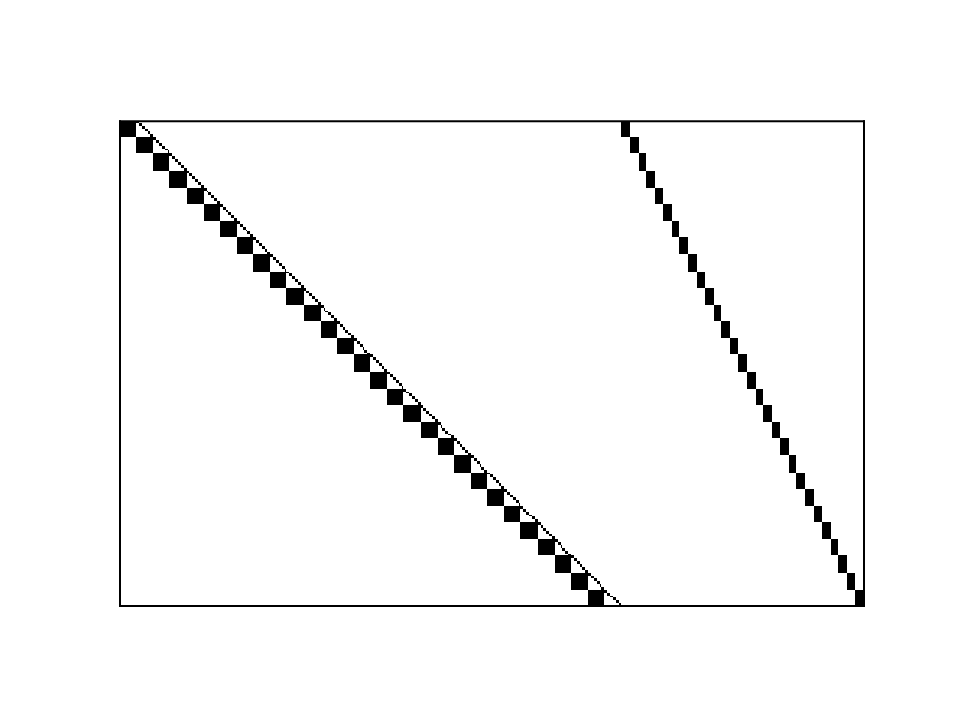}

    \vspace{-1cm}
    
    \caption{Sparsity plot of $H$ for the oscillating mass problem}
    \label{fig:H}
\end{figure}

\begin{figure}[!htb]
    \centering
    \includegraphics[width=\linewidth]{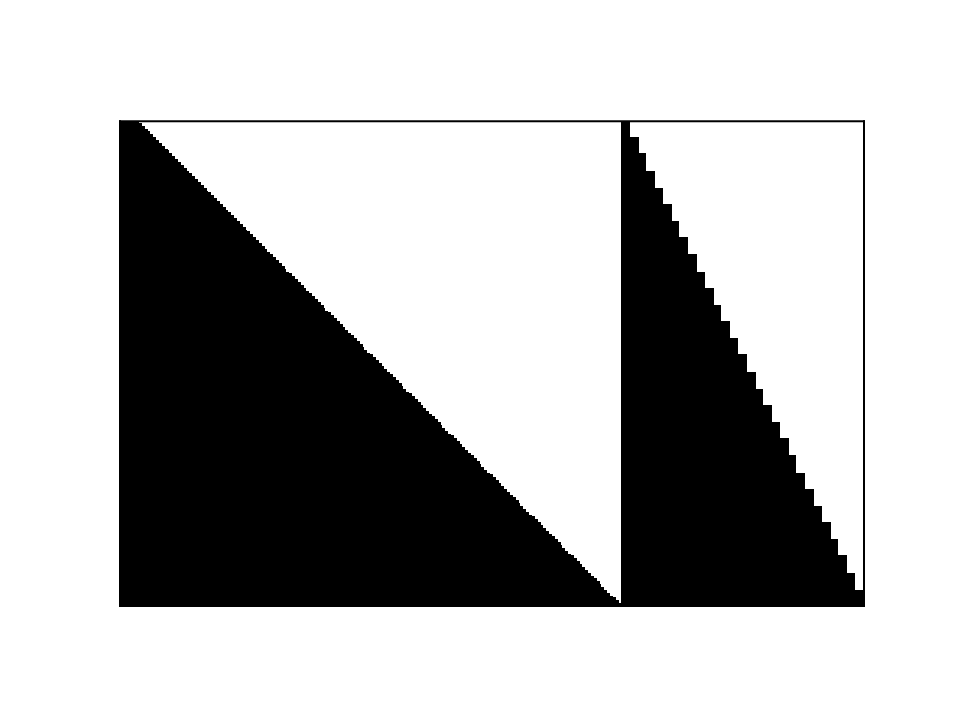}

    \vspace{-1cm}
    
    \caption{Sparsity plot of $\hat{H}$ for the oscillating mass problem}
    \label{fig:Hhat}
\end{figure}

Figures \ref{fig:H} and \ref{fig:Hhat} depict the sparsity pattern of $H$ for MPC problems before and after applying the QR preconditioner. Applying the preconditioner results in significant fill-in, i.e., the introduction of nonzero elements in a sparse matrix where they did not exist previously. In the {\pipg} implementation, $H$ is stored as a sparse matrix when the QR preconditioner is not applied, but $H$ is stored as a dense matrix when the preconditioner is applied. This choice of data structure results in faster matrix-vector multiplication for each case. Without preconditioning, the percentage of nonzero elements in $H$ is small enough that sparse matrix vector multiplication is faster than dense matrix vector multiplication. However, after applying the QR preconditioner, $H$ contains too high of a percentage of nonzero elements for sparse matrix-vector multiplication to be faster than dense matrix-vector multiplication. Consequently, the QR preconditioner is limited to problems and embedded systems where the extra storage requirement for the preconditioned, dense $H$ matrix is not restrictive. Before preconditioning, when we store $H$ for both problems in compressed sparse column (CSC), it takes roughly 187 kilobytes for the oscillating masses problem and 11.6 kilobytes for the quadrotor problem. After applying our preconditioner and storing the preconditioned matrix as a dense matrix, it takes roughly 2581 kilobytes for the oscillating masses problem and 363 kilobytes of storage for the quadrotor problem. As the problem size increases, the storage requirements for the dense matrix become more demanding.

\vspace{0.5cm}

\subsection{Analysis of results}

We assess the effect of the our algorithms using two metrics: relative distance to optimum and feasibility as defined below
\begin{equation*}
    \mathrm{error}_{\mathrm{opt}} = \frac{\|z^k - z^\star\|_\infty}{\|z^\star\|_\infty}, \; \mathrm{error}_{\mathrm{feas}} = \frac{\|Hz^k - g\|_\infty}{\|z^\star\|_\infty}
\end{equation*}

\noindent where $z^\star$ is the solution to \ref{eq:cannonical-form} computed by Gurobi \cite{gurobi} using YALMIP \cite{lofberg2004yalmip}. We terminate {\pipg} when the relative error, $\mathrm{error}_{\mathrm{opt}}$, decreases below $1.0 \times 10^{-4}$, indicating the solution has a relative error of $0.01\%$.

For the oscillating masses problem, we test the QR preconditioner against modified Ruiz equilibration. For this problem, modified Ruiz equilibration preserves closed-form projection onto set $\mb{D}$ since set $\mb{D}$ is a box before and after transformation. Additionally, we test the hypersphere preconditioner (HS) on the quadrotor problem, since in this case, the transformed set $\mb{D}$ still has a closed-form projection. If our state cost matrix $Q_t$ penalized each component of velocity by different weights or the control cost matrix $R_t$ penalized each component of thrust with different weights, the HS preconditioner will induce a change of variables which will result in the transformed set $\mb{D}$ containing ellipsoids and no longer having closed-forms projections. The hypersphere preconditioner involves multiplying the preconditioned cost function by some $\lambda > 0$ similar to how the QR preconditioner multiplies the equality constraints by some $\eta > 0$, however the authors do not provide a way of selecting this $\lambda$ \cite{abhi2023customized}. In our numerical results, we choose the $\lambda$ for the hypersphere preconditioner such that the condition number of the KKT matrix given by Equation \ref{eq:cn-bound} is minimized. 



\begin{figure}[!htb]
    \centering
    \includegraphics[width=8cm]{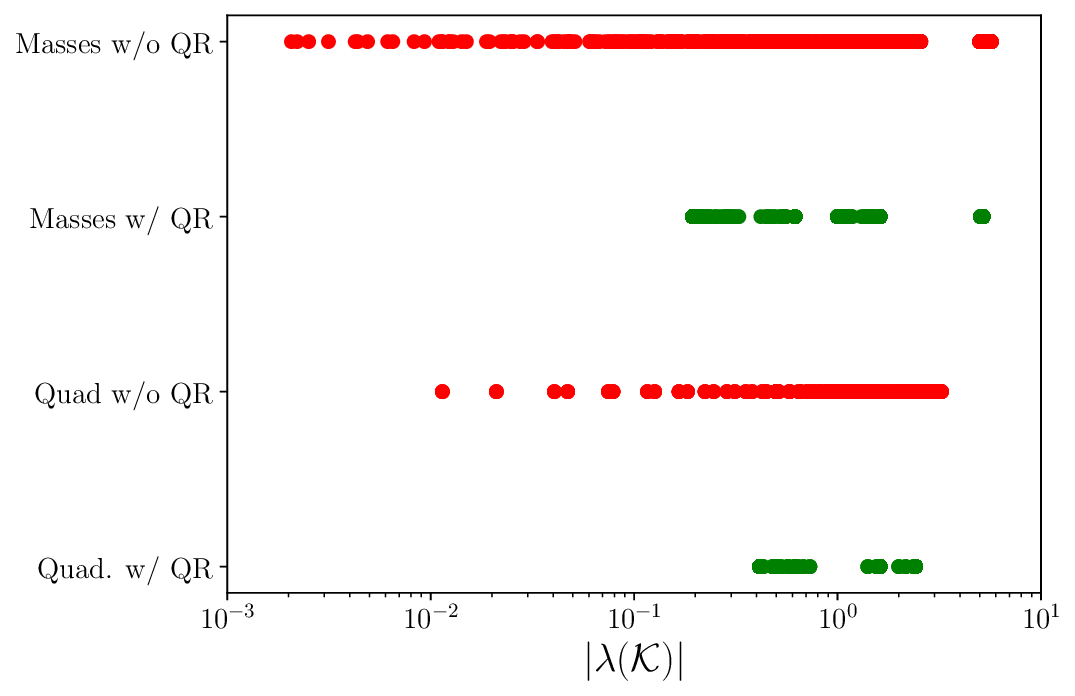}
  
    \caption{Eigenvalue distribution of the KKT matrix for the oscillating masses and quadrotor problems with and without applying the QR preconditioner}
    \label{fig:eig}
\end{figure}

Figure \ref{fig:eig} depicts the distribution of the absolute value of the eigenvalues for the KKT matrix of the oscillating masses problem and quadrotor problem before and after applying the QR preconditioner. Applying the preconditioner shrinks the interval on which the eigenvalues are distributed, making the maximum and minimum eigenvalues closer in magnitude to each other, reducing the condition number of the KKT matrix. 

\begin{figure}[!htb]
    \centering
    \includegraphics[width=8cm]{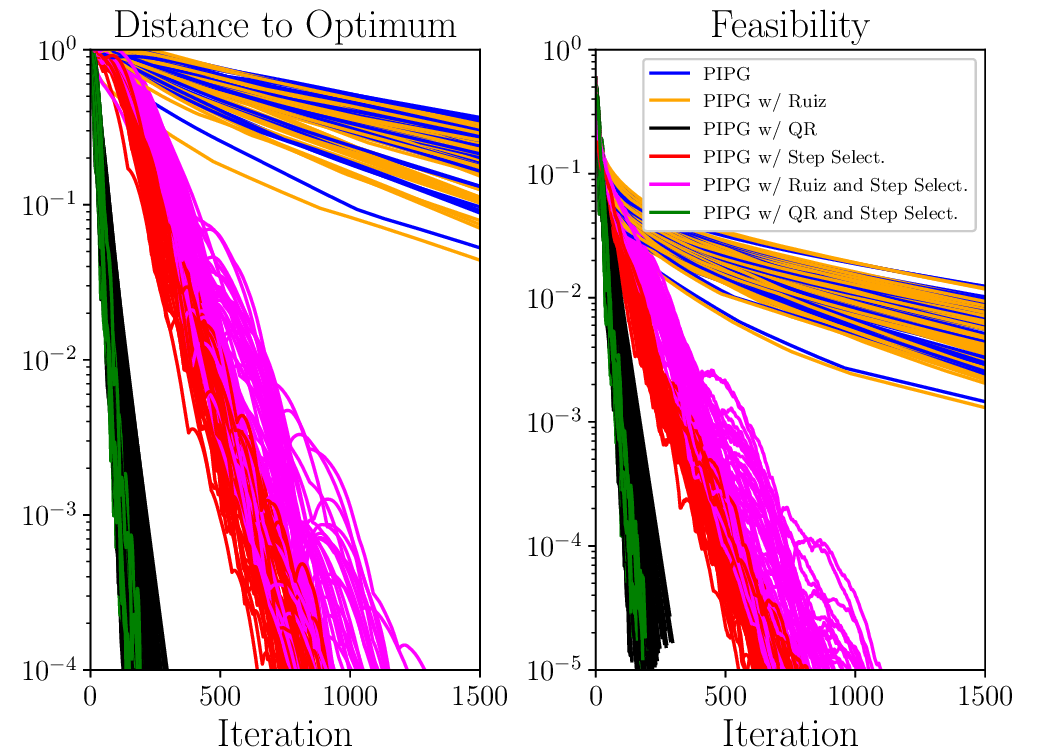}
  
    \caption{$\mathrm{error}_{\mathrm{opt}}$ (left) and $\mathrm{error}_{\mathrm{feas}}$ (right) for the oscillating masses problem}
    \label{fig:oscillating-mass-results}
\end{figure}

\begin{figure}[!htb]
    \centering
    \includegraphics[width=8cm]{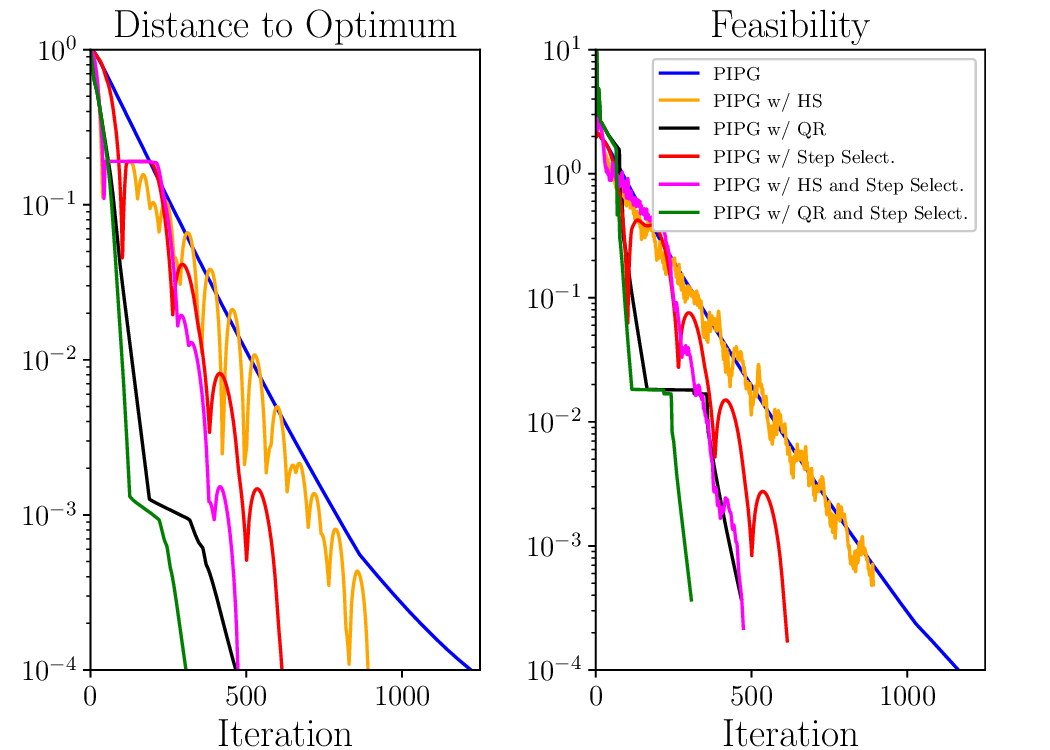}
  
    \caption{$\mathrm{error}_{\mathrm{opt}}$ (left) and $\mathrm{error}_{\mathrm{feas}}$ (right) for the quadrotor problem}
    \label{fig:quadrotor-results}
\end{figure}

Figures \ref{fig:oscillating-mass-results} and \ref{fig:quadrotor-results} show that the QR preconditioner and the step-size heuristic cause {\pipg} to converge in fewer iterations when applied independently. When applied together, we observed a further reduction in iteration count. Applying the QR preconditioner also results in fewer iterations to convergence than applying modified Ruiz equilibration or the hypersphere preconditioner. Note that for the oscillating-mass problem, {\pipg} without step-selection and the QR preconditioner requires many more than $1500$ iterations to meet our stopping criteria, but we only plotted the first $1500$ iterations for clarity.

Tables \ref{table:oscillating-mass} and \ref{table:quadcopter} contain the solve-times for our numerical examples. The solve-times in the table do not include preconditioning time for any of the preconditioners, since for MPC problems we can perform the preconditioning offline and only have to incur cost of solving the preconditioned problem. We observe that the time to perform the QR preconditioning was very small for the problems we considered: $11.9$ and $1.75$ milliseconds for the oscillating masses problem and quadrotor problem respectively. Our two proposed algorithms result in solve-time reductions, and we see that the QR preconditioner is more effective at reducing solve-times than modified Ruiz equilibration and the hypersphere preconditioner. 

\begin{table}[!htb]
\centering
\begin{tabular}{|l|l|l|l|}
\hline
               & Without Precond. & With Ruiz & With QR \\ \hline
Without Step Select. & 651.73              & 562.80           & 23.96       \\ \hline
With Step Select.    & 50.85              & 60.71           & \textbf{18.69}       \\ \hline
\end{tabular}
\caption{Average solve-time in milliseconds over 50 runs for the oscillating masses problem}
\label{table:oscillating-mass}
\end{table}

\begin{table}[!htb]
\centering
\begin{tabular}{|l|l|l|l|}
\hline
               & Without Precond. & With HS & With QR \\ \hline
Without Step Select. & 226.92              & 169.50           & 93.68       \\ \hline
With Step Select.    & 114.91              & 91.80           &  \textbf{62.04}      \\ \hline
\end{tabular}

\caption{Average solve-time in milliseconds over 100 runs for the quadcopter problem}
\label{table:quadcopter}
\end{table}

\FloatBarrier
\section{CONCLUSION}
We develop the \textit{QR preconditioner} to speed up {\pipg}, a first-order primal-dual algorithm, by applying a transformation to equality constraints which minimizes the condition number of the KKT matrix for the optimization problem. We then introduce a step-size heuristic for {\pipg} based on minimizing the upper bound on the primal-dual gap. The effectiveness of these algorithms are demonstrated on two MPC problems, where we provide comparisons to existing preconditioning algorithms.

\section*{Acknowledgments}
The authors would like to thank Dan Calderone for early conversations on the QR preconditioner, Purnanand Elango for various discussions and feedback on drafts of this paper, and Samet Uzun for his feedback on drafts.

\bibliographystyle{unsrt}
\bibliography{references}

\begin{thebibliography}{10}

\bibitem{Bauschke2011ConvexAA}
Heinz~H. Bauschke and Patrick~L. Combettes.
\newblock Convex analysis and monotone operator theory in hilbert spaces.
\newblock In {\em CMS Books in Mathematics}, 2011.

\bibitem{bauschke2018projecting}
Heinz~H. Bauschke, Minh~N. Bui, and Xianfu Wang.
\newblock Projecting onto the intersection of a cone and a sphere.
\newblock {\em SIAM Journal on Optimization}, 28(3):2158--2188, 2018.

\bibitem{jerez2014embedded}
Juan~L. Jerez, Paul~J. Goulart, Stefan Richter, George~A. Constantinides, Eric~C. Kerrigan, and Manfred Morari.
\newblock Embedded online optimization for model predictive control at megahertz rates.
\newblock {\em IEEE Transactions on Automatic Control}, 59(12):3238--3251, 2014.

\bibitem{banjac2017embedded}
Goran Banjac, Bartolomeo Stellato, Nicholas Moehle, Paul Goulart, Alberto Bemporad, and Stephen Boyd.
\newblock Embedded code generation using the osqp solver.
\newblock In {\em 2017 IEEE 56th Annual Conference on Decision and Control (CDC)}, pages 1906--1911, 2017.

\bibitem{ghadimi2015optimal}
Euhanna Ghadimi, André Teixeira, Iman Shames, and Mikael Johansson.
\newblock Optimal parameter selection for the alternating direction method of multipliers (admm): Quadratic problems.
\newblock {\em IEEE Transactions on Automatic Control}, 60(3):644--658, 2015.

\bibitem{osqp}
B.~Stellato, G.~Banjac, P.~Goulart, A.~Bemporad, and S.~Boyd.
\newblock {OSQP}: an operator splitting solver for quadratic programs.
\newblock {\em Mathematical Programming Computation}, 12(4):637--672, 2020.

\bibitem{abhi2023customized}
Abhinav~G Kamath, Purnanand Elango, Taewan Kim, Skye Mceowen, Yue Yu, John~M Carson~III, Mehran Mesbahi, and Beh{\c{c}}et A{\c{c}}{\i}kme{\c{s}}e.
\newblock Customized real-time first-order methods for onboard dual quaternion-based 6-dof powered-descent guidance.
\newblock In {\em AIAA SciTech 2023 Forum}, 2023.

\bibitem{yu2022conic}
Yue Yu, Purnanand Elango, Ufuk Topcu, and Behçet Açıkmeşe.
\newblock Proportional–integral projected gradient method for conic optimization.
\newblock {\em Automatica}, 142:110359, 2022.

\bibitem{yu2021mpc}
Yue Yu, Purnanand Elango, and Behçet Açıkmeşe.
\newblock Proportional-integral projected gradient method for model predictive control.
\newblock {\em IEEE Control Systems Letters}, 5(6):2174--2179, 2021.

\bibitem{yu2022extrapolated}
Yue Yu, Purnanand Elango, Behçet Açıkmeşe, and Ufuk Topcu.
\newblock Extrapolated proportional-integral projected gradient method for conic optimization.
\newblock {\em IEEE Control Systems Letters}, 7:73--78, 2023.

\bibitem{elango2022customized}
Purnanand Elango, Abhinav~G Kamath, Yue Yu, John~M Carson~III, Mehran Mesbahi, and Beh{\c{c}}et A{\c{c}}{\i}kme{\c{s}}e.
\newblock A customized first-order solver for real-time powered-descent guidance.
\newblock In {\em AIAA SciTech 2022 Forum}, page 0951, 2022.

\bibitem{kamath2023seco}
Abhinav~G. Kamath, Purnanand Elango, Yue Yu, Skye Mceowen, Govind~M. Chari, John~M. {Carson III}, and {Beh\c{c}et} {A\c{c}\i{}kme\c{s}e}.
\newblock Real-time sequential conic optimization for multi-phase rocket landing guidance.
\newblock {\em IFAC-PapersOnLine}, 56(2):3118--3125, 2023.
\newblock 22nd IFAC World Congress.

\bibitem{chari2024fast}
Govind~M. Chari, Abhinav~G. Kamath, Purnanand Elango, and Behcet Acikmese.
\newblock Fast monte carlo analysis for 6-dof powered-descent guidance via gpu-accelerated sequential convex programming.
\newblock In {\em AIAA SciTech 2024 Forum}, 2024.

\bibitem{chari2024spacecraft}
Govind~M Chari and Beh{\c{c}}et A{\c{c}}{\i}kme{\c{s}}e.
\newblock Spacecraft rendezvous guidance via factorization-free sequential convex programming using a first-order method.
\newblock {\em arXiv preprint arXiv:2402.04561}, 2024.

\bibitem{benzi2005numerical}
Michele Benzi, Gene~H. Golub, and Jörg Liesen.
\newblock Numerical solution of saddle point problems.
\newblock {\em Acta Numerica}, 14:1–137, 2005.

\bibitem{golub1996matrix}
Gene~H Golub and Charles~F Van~Loan.
\newblock {\em Matrix Computations}.
\newblock Johns Hopkins Series in the Mathematical Sciences. Johns Hopkins University Press, Baltimore, MD, 3 edition, October 1996.

\bibitem{nocedal2006numerical}
Jorge Nocedal and Stephen Wright.
\newblock {\em Numerical Optimization}.
\newblock Springer Series in Operations Research and Financial Engineering. Springer, New York, NY, 2 edition, July 2006.

\bibitem{wang2010fast}
Yang Wang and Stephen Boyd.
\newblock Fast model predictive control using online optimization.
\newblock {\em IEEE Transactions on Control Systems Technology}, 18(2):267--278, 2010.

\bibitem{gurobi}
{Gurobi Optimization, LLC}.
\newblock {Gurobi Optimizer Reference Manual}, 2023.

\bibitem{lofberg2004yalmip}
J.~L{\"{o}}fberg.
\newblock Yalmip : A toolbox for modeling and optimization in matlab.
\newblock In {\em In Proceedings of the CACSD Conference}, Taipei, Taiwan, 2004.

\end{thebibliography}

\end{document}